\documentclass[12pt]{amsart}

\usepackage{mathtools}
\mathtoolsset{showonlyrefs=true}
\usepackage[hmargin=0.8in,height=8.6in]{geometry}
\usepackage{amssymb,amsthm, times}
\usepackage{delarray,verbatim}
\usepackage{ifpdf}
\ifpdf
\usepackage[pdftex]{graphicx}
 \else
\usepackage[dvips]{graphicx}
 \fi

\usepackage{bm}

\usepackage{ifpdf}
\usepackage{color}
\definecolor{webgreen}{rgb}{0,.5,0}
\definecolor{webbrown}{rgb}{.8,0,0}
\definecolor{emphcolor}{rgb}{0.5,0.95,0.95}

\usepackage{hyperref}
\hypersetup{%
	colorlinks=true,
	linkcolor=webbrown,
	filecolor=webbrown,
	citecolor=webgreen,
	breaklinks=true}
\ifpdf \hypersetup{pdftex,
	pdfstartview=FitH, 
	bookmarksopen=true,
	bookmarksnumbered=true
} \else \hypersetup{dvips} \fi
\allowdisplaybreaks

\numberwithin{equation}{section}

\newtheorem{theorem}{Theorem}[section]

\newtheorem{corollary}{Corollary}[section]
\newtheorem{remark}{Remark}[section]
\newtheorem{lemma}{Lemma}[section]

\numberwithin{remark}{section} \numberwithin{proposition}{section}
\numberwithin{corollary}{section}

\newcommand{\ARXIV}[1]{\href{https://arXiv.org/abs/#1}{arXiv:#1}}

\begin{document}
	
	\title{Alpha-stable branching and beta-Frequency processes, \\beyond the IID assumption }
	\author[Adri\'an Gonz\'alez]{ Adri\'an~Gonz\'alez~Casanova$^*$}
\thanks{$*$\, Instituto de Matem\'aticas, Universidad Nacional Aut\'onoma de M\'exico, M\'exico and Department of Statistics, University of California, Berkeley. Email: adrian.gonzalez@im.unam.mx, gonzalez.casanova@berkeley.edu}
\author[Imanol~Nu\~nez]{Imanol~Nu\~nez$^\dagger$}
\thanks{$\dagger$\, Department of Probability and Statistics, Centro de Investigaci\'on en Matem\'aticas A.C. Calle Jalisco
	s/n. C.P. 36240, Guanajuato, Mexico. Email: imanol.nunez@cimat.mx, jluis.garmendia@cimat.mx.}
\author[J. L. P\'erez]{Jos\'e-Luis P\'erez$^\dagger$}

	
	\begin{abstract}
		In \cite{7authors}, Birkner et al. obtained necessary and sufficient conditions for the frequency between two independent and identically distributed continuous-state branching 
		processes time-changed by a functional of the total mass process to be 
		a Markov process. 
		Foucart et al. extended this result in \cite{Clemente} to continuous-state branching processes with immigration. 
		We generalize these results by dropping the 
		\textit{independent and identically distributed} assumption. Our result clarifies under which conditions a multi-type $\Lambda$-coalescent can be constructed from a multi-type branching process by a time change using the total mass. Finally, we address a problem formulated by Griffiths in \cite{Griff}, by clarifying the relation between 2-type $\alpha$-stable continuous-state branching processes and 2-type $\beta$-Fleming--Viot processes with mutation and selection.
			\ \ 
		\\
		\noindent \small{\noindent  \textbf{MSC2020 Subject Classifications}: 60J80; 60J90; 60G52. }
	\end{abstract}
	
	\maketitle
	\section{Introduction}
	Let $(X^1,X^2)$ be a pair of identically distributed and independent continuous-state branching processes. We can think of a population of two types of individuals in which $X^i$ represents the size of the population of type $i$ individuals for $i=1,2$. Under this setting we consider the total population size process  $Z_t:=X^1_t+X^2_t$, and if $Z_t>0$ the frequency process $R:=\{R_t:t\geq 0\}$ of type $1$ individuals given by
	\begin{align*}
		R_t=\frac{X^1_t}{X_t^1+X^2_t},\qquad t\geq0.
	\end{align*}
	Observe that the frequency process $R$ is not a Markov process. However, it is known that for the case where $X^1$ and $X^2$ are independent Feller diffusions, one can time-change the frequency process $R$ by a functional of the total population size process $Z$ to obtain the Wright--Fisher diffusion.
	
	In a more general setting, consider a function $\beta:\mathbb{R}_+\to \mathbb{R}_+$ and the following functional of the total size of the population size process $Z$ 
	\[
	T(t)=\int_0^t\beta(Z_s)ds,\qquad t\geq0.
	\]
	In \cite{7authors}, Birkner et al. characterized the class of continuous-state branching processes for which there exists a function $\beta:\mathbb{R}_+\to \mathbb{R}_+$, such that the time-changed frequency process $R_{T^{-1}(t)}$ has the Markov property. It turns out that the aforementioned class is that of the $\alpha$-stable continuous-state branching processes. Furthermore, in this case the time-changed frequency process $R_{T^{-1}(t)}$ is the moment dual of a Beta$(2-\alpha,\alpha)$-coalescent.
	This result was extended by Foucart and Hénard to continuous-state branching processes with immigration \cite{Clemente}. 
	
	The aim of this note is to generalize this result in two directions:
	\begin{enumerate}
		\item $X^1$ and $X^2$ will not be assumed to have the same distribution.
		\item $X^1$ and $X^2$ will not be assumed to be independent. Instead, we will assume that $X=(X^1,X^2)$ is given by a two-type continuous-state branching process with immigration.
	\end{enumerate}
	The frequency process associated with a pair of independent continuous-state branching processes with different distributions was studied in \cite{CGP}. On the other hand, we are not aware of studies in the literature on the frequency process arising from multi-type continuous-state branching processes.
	
	Our main result, Theorem 1.1, formulates necessary and sufficient conditions under which the frequency process $R$ can be time-changed by a functional of the total size of the
	population process Z to obtain a Markov process, and to characterize the associated frequency process.
	For simplicity, we work with real-valued processes, as opposed to \cite{7authors} and \cite{Clemente} where the authors deal with measure-valued processes.  
	
	Finally, we address the problem formulated by Griffiths in the concluding section of \cite{Griff}, about the relation between multi-type $\Lambda$-Fleming--Viot processes with mutation and selection and branching processes; by showing (see Remark \ref{Griffiths}) that the time-changed frequency process associated to a 2-type $\alpha$-stable continuous-state branching process with $\alpha\in(0,1)$ corresponds to the 2-type $\beta$-Fleming--Viot process with mutation and selection. 
	
	\subsection{Multi-type continuous-state branching processes with immigration}
	
	Denote by $D:=\{(x_1,x_2)\in\mathbb{R}^2: x_1,x_2\geq 0\}$ and $\xi_1,\xi_2\in\mathcal{C}_0^{\infty}(\mathbb{R}^2)$ such that $\xi_i=x_i$ for $i=1,2$ on some neighborhood $U$ of the origin. Under this setting, we consider a two-type continuous-state branching process with immigration $X=\{(X^1_t,X^2_t):t\geq 0\}$.  Define 
	\begin{align*}
		\mathcal{C}_0(D):&=\{f(x):\text{continuous on $D$ such that $\lim_{|x|\to\infty}f(x)=0$}\},\notag\\
		\mathcal{C}_0^n(D):&=\{f(x) \in C_0(D):\text{all the derivatives up to $n$-th order are in $\mathcal{C}_0(D)$}\},
	\end{align*}
	and $\mathcal{C}_0^{\infty}(D):=\cap_{n=1}^{\infty}\mathcal{C}_0^n(D)$.
	
	Following \cite{W} we have that $X$ is a Feller process with infinitesimal generator $\mathcal{L}$ given, for any $f\in \mathcal{C}_0^2(D)$, by
	\begin{align}\label{inf_gen_MCBI}
		&\mathcal{L}f(x)=c_1x_1\partial_{11}f(x)+c_2x_2\partial_{22}f(x)+(b_{11}x_1+b_{21}x_2+\eta_1)\partial_1f(x)+(b_{22}x_2+b_{12}x_1+\eta_2)\partial_2f(x)\notag\\
		&+\sum_{i=1}^2x_i\int_{D}\left[f(x+y)-f(x)-\xi_i(y)\partial_if(x)\right]m^i(dy)+\int_{D}\left[f(x+y)-f(x)\right]\nu(dy),
	\end{align}
	where, $c_1,c_2,b_{12},b_{21},\eta_1,\eta_2\in\mathbb{R}_+$, $b_{11},b_{22}\in\mathbb{R}$, 
	\begin{align*}
		\int_{U}&\left(\xi_1^2(y)+\xi_2(y)\right)m^1(dy)+\int_{U}\left(\xi_1(y)+\xi_2^2(y)\right)m^2(dy)\notag\\&+\sum_{i=1}^2m^i(D\backslash U)+\int_{U}\left(\xi_1(y)+\xi_2(y)\right)\nu(dy)+\nu(U\backslash U)<\infty,
	\end{align*}
	and for $i=1,2$, $\xi_i\in \mathcal{C}_0^{\infty}(D)$ such that $\xi_i=x_i$ in a neighborhood $U$ of the origin $O$.
	\subsection{Main result}
	We consider a population of two types of individuals $1$ and $2$, and the dynamics of the population of type $i$ is given by the process $X^i$ for each $i=1,2$. In order to study the evolution of the population, we define the total size of the population process
	\[
	Z_t:=X^1_t+X^2_t,\quad t\geq0, \ Z_0=z,
	\] 
	where $z:=x_1+x_2$.
	Our goal is to study the frequency process $R:=\{R_t:t\geq 0\}$ of type $1$ individuals given by
	\begin{align*}
		R_t=\frac{X^1_t}{X_t^1+X^2_t}1_{\{t<\tau\}}+\Delta1_{\{t\geq \tau\}},\quad t\geq0, \ R_0=r,
	\end{align*}
	where $r:=x_1/(x_1+x_2)$, 
	\[
	\tau:=\inf\{t\geq0:Z_t=0\}\wedge\inf\{t\geq0:Z_t=\infty\},
	\]
	and $\Delta$ is a cemetery state.
	We consider a mapping $\beta:\mathbb{R}_+\to\mathbb{R}_+$ and define a functional $T$ of the total size population process $Z$ by
	\begin{equation}\label{beta}
		T(t)=\int_0^t\beta(Z_s)ds,\qquad t\geq0.
	\end{equation}
	Additionally, we denote by $T^{-1}$ the right-continuous inverse of $T$, i.e. $T^{-1}(t):=\inf\{s:T(s)>t\}$.
	
	We are interested in studying the evolution of the frequency process $R$, however, although the pair $(R, Z)$ is a Markov process, the process $R$ is not markovian.
	Therefore in the spirit of Theorem 1.1 in \cite{7authors}, we use the random time-change $T^{-1}$ to obtain an autonomous Markov process that describes the frequency of type 1 individuals in the population. 
	This is stated in the following result, which generalizes the random-time change technique from \cite{7authors} to the setting where the processes $X^1$ and $X^2$ are dependent, not identically distributed, and where immigration is included in the model.
	
	Through the rest of the paper we denote by $\mathbb{S}^2$ to the unit sphere in $\mathbb{R}^2$ and $\mathbb{S}^2_+:=\{\xi=(\xi_1,\xi_2)\in \mathbb{S}^2: \xi_1,\xi_2\geq0\}$.
	\begin{theorem}\label{main}
		The frequency process $R$ can be time-changed by a functional of the total size of the population process $Z$ to obtain a Markov process if and only if 
		\begin{itemize}
			\item[(i)] \textbf{Continuous case.-} $m^i=\nu=b_{11}=b_{12}=b_{21}=b_{22}=0$,
			\item[(ii)] \textbf{Independent branching and multi-type immigration.-} 
			$c_1=c_2=\eta_1=\eta_2=b_{21}=b_{12}=0$,
			additionally
			\begin{align*}
				b_{ii}=\int_D\left(\xi_i(u)-u_i\right)m^i(du), \qquad \text{$i=1,2$,}
			\end{align*}
			where for $A\in\mathcal{B}(\mathbb{R}^2_+)$
			\begin{align*}
				\nu(A)=\int_{\mathbb{S}^2_+}\lambda^I(d\xi)\int_0^{\infty}1_A(r\xi) \frac{dr}{r^{\alpha}},\qquad m^i(A)=a_i\int_0^{\infty}1_A(re_i)\frac{dr}{r^{1+\alpha}}, \qquad \text{$i=1,2$,}
			\end{align*}
			with $\alpha\in(1,2)$, $e_1=(1,0)$, $e_2=(0,1)$, $\lambda^I$ is finite measure on $\mathbb{S}^2_+$, and $a_i\geq 0$ for $i=1,2$.
			\item[(iii)] \textbf{Multi-type branching.-} $\nu = c_1=c_2=\eta_1=\eta_2=b_{21}=b_{12}=0$,
			\begin{align*}
				b_{ii}=\int_{D}\xi_i(u)m^i(du), \qquad 
				m^i(A)=\int_{\mathbb{S}^2_+}\lambda^i(d\xi)\int_0^{\infty}1_A(r\xi)\frac{dr}{r^{1+\alpha}},\qquad i=1,2, 
			\end{align*}
			with $\alpha\in(0,1)$, and $\lambda^i$ is a finite measure on $\mathbb{S}^2_+$ for $i=1,2$.
		\end{itemize}
	\end{theorem}
	\begin{remark}
		Theorem \ref{main} explores different cases in which the frequency process associated with a pair of continuous-state branching processes can be transformed into a Markov process by a time-change involving a functional of the total mass. In each case, the specific conditions for this transformation depend on the particular characteristics of the branching processes involved.
		
		In case (i), the theorem considers continuous-state branching processes with immigration that are both continuous and independent, and may have different distributions. The condition for transforming the frequency process into a Markov process depends on the properties of the immigration process and the branching mechanisms involved.
		
		In case (ii), the theorem looks at independent continuous-state branching processes with multi-type immigration. The branching mechanisms in this case must be associated with stable Lévy measures of the same index for the transformation to be possible.
		
		Finally, in the last case, the theorem focuses on the frequency of one of the types in a two-type continuous-state branching process. The condition for transforming the frequency process into a Markov process involves both a stable Lévy measure governing the size of the reproduction events and a measure on the sphere that determines the amount of mass that each type receives at each reproduction event.
		
		Overall, Theorem 1.1 provides a framework for understanding the conditions under which the frequency process associated with continuous-state branching processes can be transformed into a Markov process, highlighting the role played by immigration, branching mechanisms, and other factors in this transformation.
	\end{remark}
	The proof of Theorem \ref{main} relies on the following result which is a generalization of Lemma 3.5 in \cite{7authors} to the two-dimensional case, its proof is deferred to Section \ref{proof_measure}. 
	\begin{lemma}\label{measure}
		Let $\nu$ be a measure on $D$ satisfying 
		\begin{equation}\label{integrability}
			\int_U\left(\xi_1^2(y)+\xi_2(y)\right)\nu(dy)+\nu(D\backslash U)<\infty.
		\end{equation}
		For $z>0$ let $\mu_z=\phi_z(\nu)$ be the image of $\nu$ under the mapping given by
		\begin{equation}\label{mapping}
			\phi_z:(u_1,u_2)\mapsto (r_1,r_2):=\left(\frac{u_1}{z+u_1+u_2},\frac{u_2}{z+u_1+u_2}\right).
		\end{equation}
		There exists a measure $\mu$ in $D\backslash\{(0,0)\}$ and a measurable mapping $\beta:\mathbb{R}_+\mapsto\mathbb{R}_+$ such that
		\begin{equation}\label{stameas_1}
			\mu_z=\beta(z)\mu
		\end{equation}
		if and only if, for some $\alpha\in(0,2)$,
		\begin{equation}\label{stameas}
			\nu(B)=\int_{\mathbb{S}^2_+}\lambda(d\xi)\int_0^{\infty}1_B(r\xi)\frac{dr}{r^{1+\alpha}}, \qquad B\in\mathcal{B}(\mathbb{R}^2),
		\end{equation}
		where $\lambda$ is a finite measure on $\mathbb{S}^2_+$ and $\beta(z)=\text{const}\cdot z^{-\alpha}$. Additionally, for the case $\alpha\in(1,2)$ we have that $\lambda(\mathbb{S}^2_+\backslash\{e_1\})=0$. 
	\end{lemma}
	We end this section by providing a characterization of the time change $T(t)$ and the resulting time-changed frequency process $\overline{R}:=\{R_{T^{-1}(t)}:t\geq0\}$ hinted by Theorem \ref{main} in each of the three cases. This is the content of the following corollaries and their proof is contained in the proof of Theorem \ref{main}.
	\begin{corollary}\label{cor_1}
		For case (i) in Theorem \ref{main} let $T(t):=\int_0^tZ_s^{-1}ds$ for $t>0$ (i.e. $\beta(z)=z^{-1}$ in \eqref{beta}). Then the process $\overline{R}$ is the unique weak solution to the following stochastic differential equation
		\begin{align}\label{sde_1}
			d\overline{R}_t&=2(c_2-c_1)\overline{R}_t(1-\overline{R}_t)dt+\left(\eta_1(1-\overline{R}_t)-\eta_2\overline{R}_t\right)dt\notag\\&+\sqrt{2c_1\overline{R}_t(1-\overline{R}_t)^2+2c_2(1-\overline{R}_t)\overline{R}_t^2}dB_t, \qquad t\geq0,
		\end{align}
		where $B=\{B_t:t\geq0\}$ is a Brownian motion.
		
		If $T(\tau)<\infty$, we consider $(R_{T^{-1}(t)})_{t\geq0}$ to be extended for $t \geq T(\tau)$ by an independent copy of the solution to \eqref{sde_1} started from $R_{T^{-1}(T(\tau)-)}$.
	\end{corollary}
	\begin{corollary}\label{cor_2}
		For case (ii) in Theorem \ref{main} let $\alpha\in(1,2)$ and $T(t):=\int_0^tZ_s^{1-\alpha}ds$ for $t>0$ (i.e. $\beta(z)=z^{1-\alpha}$ in \eqref{beta}). Then the process $\overline{R}$ is the unique weak solution to the following stochastic differential equation
		\begin{align} \label{new_sde_2}
			d\overline{R}_t
			& = \overline{R}_t (1 - \overline{R}_t) (a_2 - a_1) \int_{(0, 1)} 
			w^{1-\alpha} (1 - w)^{\alpha - 2} dw dt \notag\\
			& + \int_0^1 \int_0^1 w (1 - \overline{R}_{t-}) 1_{\{u \leq \overline{R}_{t-}\}} 
			\tilde{N}^1(dt, dw, du) 
			+ \int_0^1 \int_0^1 w \overline{R}_{t-} 1_{\{u \leq (1 - \overline{R}_{t-})\}} 
			\tilde{N}^2(dt, dw, du)  \notag \\
			& + \int_0^1 \int_{\mathbb{S}^2_+} y \left( \frac{\langle \xi, e_1 \rangle}{\langle \xi, 1 \rangle}
			- \overline{R}_{t-} \right) \bar{N}^{I}(dt, dy, d\xi), 
		\end{align}
		where for $i=1,2$, $\tilde{N}^i $ are compensated Poisson random measures on $(0,\infty)\times(0,1)\times (0,1)$ with intensity measures $dt\Lambda^i(dw)du$ respectively, with
		\begin{align*}
			\Lambda^i(dw)=a_i(1-w)^{\alpha-1}w^{-(1+\alpha)}dw,
		\end{align*} 
		whereas $\bar{N}^I$ 
		is a Poisson random measure on $(0, \infty) \times (0, 1) \times \mathbb{S}^2_+$ with 
		intensity measure $dt \bar{\Lambda}^I(dy,d \xi)$, with 
		\[
		\bar{\Lambda}^I(dy,d \xi) = (1 - y)^{\alpha - 2} y^{-\alpha} 
		\langle \xi, 1 \rangle^{\alpha - 1} \lambda^I(d \xi)dy .
		\]
		If $T(\tau)<\infty$, we consider $(R_{T^{-1}(t)})_{t\geq0}$ to be extended for $t \geq T(\tau)$ by an independent copy of the solution to \eqref{new_sde_2} started from $R_{T^{-1}(T(\tau)-)}$.
	\end{corollary}
	\begin{corollary}\label{cor_3}
		For case (iii) in Theorem \ref{main} let $\alpha\in(0,1)$ and $T(t):=\int_0^tZ_s^{1-\alpha}ds$ for $t>0$ (i.e. $\beta(z)=z^{1-\alpha}$ in \eqref{beta}) and $T^{-1}(t):=\inf\{s:T(s)>t\}$. Then the process $\overline{R}$ is the unique weak solution to the following stochastic differential equation
		\begin{align}\label{new_sde_3}
			d \overline{R}_t 
			& = \int_0^1 \int_{\mathbb{S}^2_+} \int_0^1 y \left( 
			\frac{\langle \xi, e_1 \rangle}{\langle \xi, 1 \rangle} - \overline{R}_{t-} 
			\right) 1_{\{u \leq \overline{R}_{t-}\}} \mathring{N}^1(dt, dy, d\xi, du) \notag\\
			& + \int_0^1 \int_{\mathbb{S}^2_+} \int_0^1 y \left( 
			\frac{\langle \xi, e_1 \rangle}{\langle \xi, 1 \rangle} - \overline{R}_{t-} 
			\right) 1_{\{u \leq (1-\overline{R}_{t-})\}} \mathring{N}^2(dt, dy, d\xi, du),
		\end{align}
		where $\mathring{N}^i$, $i = 1, 2$, are independent Poisson random measures on 
		$(0, \infty) \times (0, 1) \times \mathbb{S}^2_+ \times (0, 1)$ with intensity measures 
		$dt \bar{\Lambda}^i(dy,d\xi)du$ where 
		\[
		\bar{\Lambda}^i(dy,d\xi)
		= (1 - y)^{\alpha - 1} y^{-1 - \alpha} 
		\langle \xi, 1 \rangle^{\alpha} \lambda^i(d\xi)dy.
		\]
		If $T(\tau)<\infty$, we consider $(R_{T^{-1}(t)})_{t\geq0}$ to be extended for $t \geq T(\tau)$ by an independent copy of the solution to \eqref{new_sde_3} started from $R_{T^{-1}(T(\tau)-)}$.
	\end{corollary}
	\begin{remark}
		For some choices of the parameters, the processes on the three previous corollaries have moment duals. We now provide a few examples:
		\begin{itemize}
			\item The process $\overline{R}$ in Corollary 1.1 is known as the Gillespie--Wright--Fisher diffusion (see \cite{Gill73,Gill74}), whose dual is the branching coalescing pairwise branching process studied in \cite{GPP,GMP}. 
			\item By taking $a_1=a_2=a$,  $\lambda^I=0$ in Corollary \ref{cor_2}, and $\lambda^i(d\xi)=a\delta_{e_i}(d\xi)$ with $a\geq0$ in Corollary \ref{cor_3}, the moment dual of the frequency process $\overline{R}$ is given by the $\beta$-coalescent studied in \cite{7authors}. If additionally we take $\lambda^I(d\xi)=b_1\delta_{e_1}(d\xi)+b_2\delta_{e_2}(d\xi)$ in Corollary \ref{cor_2} we obtain the $\beta$-coalescent with coordinated mutations (see \cite{Clemente}).
			\item We want to point out that some of the moment duals that emerge seem to be new in the literature. For instance, by taking $a_1>a_2$ in Corollary \ref{cor_2}, the moment dual of the process $\overline{R}$ is given by a $\beta$-coalescent with selection and coordinated multi-type mutations, where both types participate in a mutation event. On the other hand, by taking $\lambda_1=\lambda_2$ in Corollary \ref{cor_3}, the moment dual becomes a pure mutation process. Finally, the moment dual of the frequency process given in Corollary \ref{cor_3} is related to the multi-type $\Lambda$-coalescent with selection and mutation given in Theorem 2.2 in \cite{Griff}.
		\end{itemize}
	\end{remark}
	\section{Proofs the main results.}
	\subsection{Proof of Lemma \ref{measure}}\label{proof_measure}
	Proceeding like in the proof of Lemma 3.5 in \cite{7authors}, for $c>0$ consider the mapping $\psi_c:(u_1,u_2)\mapsto c(u_1,u_2)$. Hence, for all $z,c>0$
	\begin{align}\label{stameas_2}
		\phi_z=\phi_{cz}\circ\psi_c
	\end{align}
	Using \eqref{stameas_1} together with \eqref{stameas_2} 
	\begin{align}\label{stameas_4}
		\beta(z)\mu=\phi_{cz}(\psi_c(\nu)).
	\end{align}
	Additionally, by \eqref{stameas_1}
	\begin{align}\label{stameas_3}
		\phi_{cz}(\nu)=\beta(cz)\mu
	\end{align}
	Inverting \eqref{stameas_3} we obtain
	\begin{align}\label{stameas_5}
		\nu=\beta(cz)\phi^{-1}_{cz}(\mu).
	\end{align}
	Inverting \eqref{stameas_4} and using \eqref{stameas_5}
	\begin{align}\label{stameas_6}
		\psi_c(\nu)=\frac{\beta(z)}{\beta(cz)}\nu.
	\end{align}
	Let us consider $C\in\mathcal{B}(\mathbb{S}^2_+)$, then by taking $z=1$ in \eqref{stameas_6} we have for $h>0$
	\begin{align*}
		\nu\left(\frac{1}{c}(h,\infty)\times C\right)=\psi_c(\nu)\left((h,\infty)\times C\right)=\tilde{\beta}\left(\frac{1}{c}\right)\nu\left((h,\infty)\times C\right),
	\end{align*}
	where $\tilde{\beta}\left(\frac{1}{c}\right)=\beta(1)/\beta(c)$. Then proceeding like in the proof of Lemma 3.5 in \cite{7authors} we have that
	\begin{align}\label{stameas_7}
		\nu\left((h,\infty)\times C\right)=\text{const}\cdot h^{-\alpha}.
	\end{align}
	Now let us define $\lambda (C):=\alpha\nu\left((1,\infty)\times C\right)$, then using \eqref{stameas_7}
	\begin{align*}
		\nu\left((h,\infty)\times C\right)=\frac{\lambda(C)}{\alpha} h^{-\alpha}.
	\end{align*}
	Let us denote by $\nu'$ to the right-hand side of \eqref{stameas} then for any $C\in\mathcal{B}(\mathbb{S}^2_+)$
	\begin{align*}
		\nu'((h,\infty)\times C))=\lambda(C)\int_{h}^{\infty}\frac{dr}{r^{1+\alpha}}=h^{-\alpha}\frac{\lambda(C)}{\alpha}=\nu((h,\infty)\times C).
	\end{align*}
	Hence, by an application of Dynkin's Lemma (see Theorem 3.2 in \cite{bi}) we have that $\nu=\nu'$ on $\mathcal{B}(\mathbb{R}^2\backslash\{(0,0)\})$. When $\alpha\in(1,2)$ the fact that $\lambda(\mathbb{S}^2_+\backslash\{e_1\})=0$ follows from \eqref{integrability}.
	
	\subsection{Proof of Theorem \ref{main}}
	Let us consider 
	$
	\tau_M:=\inf\{t>0: Z_t>M\}\wedge\inf\{t>0: Z_t<1/M\},
	$
	the minimum between the first hitting time of $1/M$ and the first passage time above $M$ for the process $Z$.
	Additionally, we denote by $(\mathcal{F}_t)_{t\geq0}$ the filtration generated by the process $X$.
	
	Then, using Dynkin's formula (see  Proposition IV.1.7 in \cite{EK}) we have that, for any $f\in \mathcal{C}^2([0,1])$, the process
	\[
	f(R_{t\wedge\tau_M})-f(R_0)-\int_0^{t\wedge\tau_M}\mathcal{L}f(g(X_s))ds, \qquad t\geq0,
	\]
	is a $\mathcal{F}_{t\wedge\tau}$-martingale, where 
	$
	g(x)=\displaystyle\frac{x_1}{x_1+x_2}
	$
	for $x\in D$. 
	
	Using \eqref{inf_gen_MCBI} we can write
	\begin{align}\label{inf_gen_r_1}
		&\mathcal{L}f(g(x))=c_1x_1\left[f''\left(\frac{x_1}{x_1+x_2}\right)\frac{x_2^2}{(x_1+x_2)^4}-2f'\left(\frac{x_1}{x_1+x_2}\right)\frac{x_2}{(x_1+x_2)^3}\right]\notag\\
		&+c_2x_2\left[f''\left(\frac{x_1}{x_1+x_2}\right)\frac{x_1^2}{(x_1+x_2)^4}+2f'\left(\frac{x_1}{x_1+x_2}\right)\frac{x_1}{(x_1+x_2)^3}\right]\notag\\
		&+f'\left(\frac{x_1}{x_1+x_2}\right)\frac{x_2}{(x_1+x_2)^2}\left(b_{11}x_1+b_{21}x_2+\eta_1\right)\notag\\
		&-f'\left(\frac{x_1}{x_1+x_2}\right)\frac{x_1}{(x_1+x_2)^2}\left(b_{12}x_1+b_{22}x_2+\eta_2\right)\notag\\
		&+x_1\int_D\left[f\left(\frac{x_1+u_1}{x_1+x_2+u_1+u_2}\right)-f\left(\frac{x_1}{x_1+x_2}\right)-\xi_1(u)f'\left(\frac{x_1}{x_1+x_2}\right)\frac{x_2}{(x_1+x_2)^2}\right]m^1(du)\notag\\
		&+x_2\int_D\left[f\left(\frac{x_1+u_1}{x_1+x_2+u_1+u_2}\right)-f\left(\frac{x_1}{x_1+x_2}\right)+\xi_2(u)f'\left(\frac{x_1}{x_1+x_2}\right)\frac{x_1}{(x_1+x_2)^2}\right]m^2(du)\notag\\
		&+\int_D\left[f\left(\frac{x_1+u_1}{x_1+x_2+u_1+u_2}\right)-f\left(\frac{x_1}{x_1+x_2}\right)\right]\nu(du).
	\end{align}
	By making the change of variable
	$
	\left(\frac{x_1}{x_1+x_2},x_1+x_2\right)\mapsto (r,z),
	$
	in \eqref{inf_gen_r_1}, we obtain that 
	the process 
	\begin{align*}
		M_t^{f,M}:=f(R_{t\wedge\tau_M})-f(R_0)-\int_0^{t\wedge\tau_M}\mathcal{G}(R_s,Z_s)ds, \qquad t\geq0,
	\end{align*}
	is a$\mathcal{F}_{t\wedge\tau}$-martingale, where for $r\in[0,1]$ and $z\in(0,\infty)$,
	\begin{align*}
		\mathcal{G}(r,z):&=c_1\frac{r}{z}\left[f''(r)(1-r)^2-2f'(r)(1-r)\right]+c_2\frac{(1-r)}{z}\left[f''(r)r^2+2f'(r)r\right]\notag\\&+(1-r)f'(r)\left[b_{11}r+b_{21}(1-r)+\frac{\eta_1}{z}\right]-rf'(r)\left[b_{12}r+b_{22}(1-r)+\frac{\eta_2}{z}\right]\notag\\
		&+rz\int_\Delta\Bigg[f\left(r(1-w_1-w_2)+w_1\right)-f(r)\notag\\&\hspace{3cm}-\xi_1\left(\frac{zw_1}{1-w_1-w_2},\frac{zw_2}{1-w_1-w_2}\right)f'(r)\frac{(1-r)}{z}\Bigg]\phi_z(m^1)(dw)\notag\\
		&+(1-r)z\int_\Delta\Bigg[f\left(r(1-w_1-w_2)+w_1\right)-f(r)\notag\\&\hspace{4cm}+\xi_2\left(\frac{zw_1}{1-w_1-w_2},\frac{zw_2}{1-w_1-w_2}\right)f'(r)\frac{r}{z}\Bigg]\phi_z(m^2)(dw)\notag\\
		&+\int_\Delta\left[f\left(r(1-w_1-w_2)+w_1\right)-f(r)\right]\phi_z(\nu)(dw),
	\end{align*}
	where $\Delta = \{(w_1, w_2) \in \mathbb{R}_+^2 : w_1 + w_2 < 1\}$, and where we recall that for any measure $\zeta$ on $D$, $\phi_z(\zeta)$ is the image of $\zeta$ under the mapping given in \eqref{mapping}.
	
	Recall the time-change $T$ given in \eqref{beta}, and define $\overline{R}_t:=R_{T^{-1}(t)}$ for $t\geq0$. 
	Then, by noting that $T^{-1}(T(\tau_M)) = \tau_M$ and $T^{-1}(t \wedge T(\tau_M)) = T^{-1}(t) \wedge \tau_M$, the optional stopping theorem (see Theorem 6.29 in \cite{Ka}) implies that 
	\begin{align}\label{mart_prob_3_alt}
		\overline{M}_t^{f, M} :&= f(R_{T^{-1}(t) \wedge \tau_M}) - f(R_0) - \int_0^{T^{-1}(t) \wedge \tau_M} \mathcal{G}(R_s, Z_s) ds \notag \\
		& = f(\overline{R}_{t \wedge T(\tau_M)}) - f(\overline{R}_0) - \int_0^{t \wedge T(\tau_M)} \mathcal{G}(\overline{R}_s, Z_{T^{-1}(s)}) \frac{1}{\beta(Z_{T^{-1}(s)})} ds \qquad t \geq 0,
	\end{align}
	is a $\mathcal{F}_{T^{-1}(t\wedge(T(\tau)-))}$-martingale.
	
	Now, for $s\geq 0$
	\begin{align*}
		\mathcal{G}(\overline{R}_s,&Z_{T^{-1}(s)})\frac{1}{\beta(Z_{T^{-1}(s)})}=c_1\frac{\overline{R}_s}{Z_{T^{-1}(s)}}\left[f''(\overline{R}_s)(1-\overline{R}_s)^2-2f'(\overline{R}_s)(1-\overline{R}_s)\right]\frac{1}{\beta(Z_{T^{-1}(s)})}\notag\\
		&+c_2\frac{(1-\overline{R}_s)}{Z_{T^{-1}(s)}}\left[f''(\overline{R}_s)\overline{R}_s^2+2f'(\overline{R}_s)\overline{R}_s\right]\frac{1}{\beta(Z_{T^{-1}(s)})}\notag\\
		&+(1-\overline{R}_s)f'(\overline{R}_s)\left(b_{11}\overline{R}_s+b_{21}(1-\overline{R}_s)+\frac{\eta_1}{Z_{T^{-1}(s)}}\right)\frac{1}{\beta(Z_{T^{-1}(s)})}\notag\\
		&-\overline{R}_sf'(\overline{R}_s)\left(b_{12}\overline{R}_s+b_{22}(1-\overline{R}_s)+\frac{\eta_2}{Z_{T^{-1}(s)}}\right)\frac{1}{\beta(Z_{T^{-1}(s)})}\notag\\
		&+\overline{R}_sZ_{T^{-1}(s)}\int_\Delta\Bigg[f(\overline{R}_s(1-w_1-w_2)+w_1)-f(\overline{R}_s)\notag\\
		&-\xi_1\left(\frac{Z_{\tau^{-1}(s)}w_1}{1-w_1-w_2},\frac{Z_{\tau^{-1}(s)}w_2}{1-w_1-w_2}\right)f'(\overline{R}_s)\frac{(1-\overline{R}_s)}{Z_{T^{-1}(s)}}\Bigg]\frac{1}{\beta(Z_{T^{-1}(s)})}\phi_{Z_{T^{-1}(s)}}(m^1)(dw)\notag\\
		&+(1-\overline{R}_s)Z_{T^{-1}(s)}\int_\Delta\Bigg[f(\overline{R}_s(1-w_1-w_2)+w_1)-f(\overline{R}_s)
	\end{align*}
	\begin{align}\label{time_change_gen}
		&+\xi_2\left(\frac{Z_{T^{-1}(s)}w_1}{1-w_1-w_2},\frac{Z_{\tau^{-1}(s)}w_2}{1-w_1-w_2}\right)f'(\overline{R}_s)\frac{\overline{R}_s}{Z_{T^{-1}(s)}}\Bigg]\frac{1}{\beta(Z_{T^{-1}(s)})}\phi_{Z_{T^{-1}(s)}}(m^2)(dw)\notag\\
		&+\int_\Delta\left[f(\overline{R}_s(1-w_1-w_2)+w_1)-f(\overline{R}_s)\right]\frac{1}{\beta(Z_{T^{-1}(s)})}\phi_{Z_{T^{-1}(s)}}(\nu)(dw).
	\end{align}
	By Lemma \ref{measure}, the time-changed process $\overline{R}$ will be an autonomous frequency process (with no dependence on the process $Z$) if and only if one of the three following conditions hold:
	\begin{itemize}
		\item[(i)] $\beta(z)=z^{-1}$ for $z\geq0$, and
		\begin{align*}
			m^i(du)=\nu(du)=b_{11}=b_{12}=b_{21}=b_{22}=0.
		\end{align*}
		\item[(ii)] $\beta(z)=z^{1-\alpha}$ for $z\geq0$ with $\alpha\in(1,2)$, 
		\begin{align*}
			c_1=c_2=\eta_1=\eta_2=b_{21}=b_{12}=0,
		\end{align*}
		and
		\begin{align*}
			b_{11}=\int_D\left(\xi_1(u)-u_1\right)m^1(du), \ b_{22}=\int_D\left(\xi_2(u)-u_2\right)m^2(du),
		\end{align*}
		with
		\begin{align*}
			\nu(du)=1_{r\xi}(u)\lambda^I(d\xi) \frac{dr}{r^{\alpha}},\qquad m^i(du)=a_i1_{re_i}(u)\frac{dr}{r^{1+\alpha}},\qquad i=1,2,
		\end{align*}
		where $e_1=(1,0)$ and $e_2=(0,1)$.
		
		\item[(iii)] $\beta(z)=z^{1-\alpha}$ for $z\geq0$ with $\alpha\in(0,1)$, 
		\[
		\nu(du)=c_1=c_2=\eta_1=\eta_2=b_{21}=b_{12}=0,
		\]
		and 
		\begin{align*}
			b_{11}=\int_{D}\xi_1(u)m^1(du), \qquad b_{22}=\int_D\xi_2(u)m^2(du),
		\end{align*}
		with
		\[
		m^1(du)=1_{r\xi}(u)\frac{dr}{r^{1+\alpha}}\lambda^1(d\xi), \qquad m^2(du)=1_{r\xi}(u)\frac{dr}{r^{1+\alpha}}\lambda^2(d\xi).
		\]
	\end{itemize}
	\subsubsection{Case (i)}
	For this case, using \eqref{mart_prob_3_alt} together with \eqref{time_change_gen} gives
	\begin{align}\label{mart_prob_i_a_alt}
		\overline{M}_t^{f,M}&=f(\overline{R}_{t\wedge T(\tau_M)})-f(\overline{R}_0)-\int_0^{t\wedge T(\tau_M)}f''(\overline{R}_s)\left[c_1\overline{R}_s(1-\overline{R}_s)^2+c_2(1-\overline{R}_s)\overline{R}_s^2\right]ds\notag\\
		& -\int_0^{t\wedge T(\tau_M)}2(c_2-c_1)f'(\overline{R}_s)\overline{R}_s(1-\overline{R}_s)ds - \int_0^{t\wedge T(\tau_M)}f'(\overline{R}_s)\left(\eta_1(1-\overline{R}_s)-\eta_2\overline{R}_s\right)ds,
	\end{align}
	Now, noting that $f\in\mathcal{C}^2([0,1])$, we take $M\to\infty$ in \eqref{mart_prob_3_alt} and by the dominated convergence theorem, we obtain that the right-hand side of \eqref{mart_prob_i_a_alt}, with $t\wedge T(\tau_M)$ replaced by $t\wedge (T(\tau)-)$, is a $\mathcal{F}_{T^{-1}(t\wedge (T(\tau)-))}$-martingale.
	Let us denote by $Y:=\{Y_t:t\geq0\}$ a weak solution to \eqref{sde_1}, then by Proposition 4.2 in \cite{CGP} and Proposition IV.1.7 in \cite{EK} we have that the process
	\begin{align}\label{mart_prob_i_alt}
		f(Y_{t})&-f(Y_0)-\int_0^{t}f''(Y_s)\left[c_1Y_s(1-Y_s)^2+c_2(1-Y_s)Y_s^2\right]ds\notag\\
		& -\int_0^{t}2(c_2-c_1)f'(Y_s)Y_s(1-Y_s)ds - \int_0^{t}f'(Y_s)\left(\eta_1(1-Y_s)-\eta_2Y_s\right)ds,
	\end{align}
	is a $\mathcal{F}^Y_{t}$-martingale, where $(\mathcal{F}^Y_{t})_{t\geq0}$ is the filtration generated by $Y$.
	
	By Propositions 4.1 and 4.2 in \cite{CGP} together with Proposition 4.2 in \cite{FL} the solution to the martingale problem stated in \eqref{mart_prob_i_alt} is unique. Hence, Lemma IV.5.16 in \cite{EK} implies that the process $\overline{R}$ is the unique weak solution to \eqref{sde_1}.
	\subsubsection{Case (ii)}
	By \eqref{mart_prob_3_alt} and \eqref{time_change_gen} we obtain that
	\begin{align}\label{mart_prob_ii_alt}
		&\overline{M}_t^{f,M}=f(\overline{R}_{t\wedge T(\tau_M)})-f(\overline{R}_0)\notag\\&
		-\int_0^{t\wedge T(\tau_M)} \int_0^{1}\int_{\mathbb{S}_+^2}\left[f(\overline{R}_s(1-\langle u\xi,1\rangle)+\langle u\xi,e_1\rangle)-f(\overline{R}_s)\right]1_{\{u\xi \in \Delta\}}\frac{(1- \langle\xi,1\rangle u)^{\alpha-2}}{u^{\alpha}}\lambda^I(d\xi)duds\notag\\&-\int_0^{t\wedge T(\tau_M)} \overline{R}_s\int_0^{1}\left[f(\overline{R}_s(1-w_1)+w_1)-f(\overline{R}_s)-f'(\overline{R}_s)(1-\overline{R}_s)\frac{w_1}{1-w_1}\right]a_1\frac{(1-w_1)^{\alpha-1}}{w_1^{1+\alpha}}dw_1ds\notag\\
		&-\int_0^{t\wedge T(\tau_M)} (1-\overline{R}_s)\int_0^{1}\left[f(\overline{R}_s(1-w_2))-f(\overline{R}_s)+f'(\overline{R}_s)\overline{R}_s\frac{w_2}{1-w_2}\right]a_2\frac{(1-w_2)^{\alpha-1}}{w_2^{1+\alpha}}dw_2ds.
	\end{align}
	As in case (i), taking $M\to\infty$ we obtain that the right-hand side of \eqref{mart_prob_ii_alt}, with $t\wedge T(\tau_M)$ replaced by $t\wedge (T(\tau)-)$, is a $\mathcal{F}_{T^{-1}(t\wedge (T(\tau)-))}$-martingale. By a slight modification of Proposition 4.2 in \cite{CGP}, for any solution $Y$ to \eqref{new_sde_2}, the process  
	\begin{align}\label{mart_prob_ii_1}
		f(Y_{t}&)-f(Y_0)\notag\\&
		-\int_0^{t} \int_0^{1}\int_{\mathbb{S}_+^2}\left[f\left(Y_s(1-y) + y\frac{\langle \xi,e_1\rangle}{\langle \xi,1\rangle}\right)-f(Y_s)\right]\frac{(1- y)^{\alpha-2}}{y^{\alpha}}\langle\xi,1\rangle^{\alpha-1}\lambda^I(d\xi)dyds\notag\\&-\int_0^{t} Y_s\int_0^{1}\left[f(Y_s(1-w_1)+w_1)-f(Y_s)-f'(Y_s)(1-Y_s)\frac{w_1}{1-w_1}\right]a_1\frac{(1-w_1)^{\alpha-1}}{w_1^{1+\alpha}}dw_1ds\notag\\
		&-\int_0^{t} (1-Y_s)\int_0^{1}\left[f(Y_s(1-w_2))-f(Y_s)+f'(Y_s)Y_s\frac{w_2}{1-w_2}\right]a_2\frac{(1-w_2)^{\alpha-1}}{w_2^{1+\alpha}}dw_2ds,
	\end{align}
	is a $\mathcal{F}_{t}^Y$-martingale, where $(\mathcal{F}^Y_{t})_{t\geq0}$ is the filtration generated by $Y$. 
	Notice that the change of variable $y = u \langle \xi, 1\rangle$ in the first integral of \eqref{mart_prob_ii_1} leads to a similar expression as the first integral in \eqref{mart_prob_ii_alt}.
	
	By a slight modification of Proposition 4.1 in \cite{CGP} together with Proposition 4.2 in \cite{CGP} and Proposition 4.2 in \cite{FL}, we have the martingale problem given in \eqref{mart_prob_ii_1} has a unique solution. 
	Therefore, using Lemma IV.5.16 in \cite{EK} we have that the process $\overline{R}$ is the unique weak solution to \eqref{new_sde_2}.
	\subsubsection{Case (iii)}
	
	Using \eqref{mart_prob_3_alt} together with \eqref{time_change_gen} we obtain that
	\begin{align}\label{mart_prob_iii_alt}
		&\overline{M}_t^{f,M}=f(\overline{R}_{t\wedge T(\tau_M)})-f(\overline{R}_0)\notag\\&-\int_0^{t\wedge T(\tau_M)}\overline{R}_s \int_0^{1}\int_{\mathbb{S}_+^2}\Big[f(\overline{R}_s(1-\langle u\xi,1\rangle)+\langle u\xi,e_1\rangle)\notag\\&\hspace{6cm}-f(\overline{R}_s)\Big]1_{\{u\xi \in \Delta\}}\lambda^1(d\xi)\frac{(1-\langle\xi,1\rangle u)^{\alpha-1}}{u^{1+\alpha}}duds\notag\\
		&-\int_0^{t\wedge T(\tau_M)}(1-\overline{R}_s) \int_0^{1}\int_{\mathbb{S}_+^2}\Big[f(\overline{R}_s(1-\langle u\xi,1\rangle)+\langle u\xi,e_1\rangle)\notag\\&\hspace{6cm}-f(\overline{R}_s)\Big]1_{\{u\xi \in \Delta\}}\lambda^2(d\xi)\frac{(1-\langle\xi,1\rangle u)^{\alpha-1}}{u^{1+\alpha}}duds.
	\end{align}
	Proceeding like in the previous case, by taking $M \to \infty$ we obtain that the right hand-side of \eqref{mart_prob_iii_alt}, replacing $t\wedge T(\tau_M)$ by $t\wedge (T(\tau)-)$, is a a $\mathcal{F}_{T^{-1}(t\wedge (T(\tau)-))}$-martingale. As in case (ii), by a minor modification of Proposition 4.2 in \cite{CGP} we obtain, for any solution $Y$ to \eqref{new_sde_3}, that the process
	\begin{align}\label{mart_prob_iii_a}
		f(Y_t)&-f(Y_0)\notag\\&-\int_0^{t}Y_s \int_0^{1}\int_{\mathbb{S}_+^2}\left[f\left(Y_s(1-y)+\frac{\langle \xi,e_1\rangle}{\langle \xi,1\rangle}y\right)-f(Y_s)\right] \lambda^1(d\xi)\frac{(1-y)^{\alpha-1}}{y^{1+\alpha}}\langle \xi,1\rangle^{\alpha}dyds\notag\\
		&-\int_0^{t\wedge\tau_M}(1-Y_s) \int_0^{1}\int_{\mathbb{S}_+^2}\left[f\left(Y_s(1-y)+\frac{\langle \xi,e_1\rangle}{\langle \xi,1\rangle}y\right)-f(Y_s)\right] \lambda^2(d\xi)\frac{(1-y)^{\alpha-1}}{y^{1+\alpha}}\langle \xi,1\rangle^{\alpha}dyds.
	\end{align}
	is a a $\mathcal{F}_{t}^Y$-martingale, with $(\mathcal{F}^Y_{t})_{t\geq0}$ the filtration generated by $Y$.
	
	Using Propositions 4.1 and 4.2 in \cite{CGP} and Proposition 4.2 in \cite{FL}, we obtain the uniqueness of the solution to the martingale problem \eqref{mart_prob_iii_a}. Hence, as in the previous case, an application of Lemma IV.5.16 in \cite{EK} gives that the process $\overline{R}$ is the unique weak solution to \eqref{new_sde_2}.
	
	\begin{remark}\label{Griffiths}
		In Section 2.4 in \cite{Griff}, Griffiths introduced the 2-type $\Lambda$-Fleming--Viot process with mutation and selection, $Y=\{Y_t:t\geq0\}$, through its infinitesimal generator $\mathcal{A}$, 
		given for any $g\in\mathcal{C}^2([0,1])$ by
		\begin{align*}
			\mathcal{A}g(x)=\int_{\overline{\Delta}}\sum_{i=1}^2 x^{2-i} (1-x)^{i-1} \left(g(x(1-\langle y,1\rangle)+\langle y, e_1\rangle)-g(x)\right)\Lambda^i(dy),
		\end{align*}
		where $\overline{\Delta}:=\{y\in\mathbb{R}_+^2:\langle y,1\rangle\leq 1\}$ and $\Lambda^i$ is a measure on $\overline{\Delta}$ such that $\int_{\overline{\Delta}}\langle y,1\rangle \Lambda^i(dy)<\infty$ for $i=1,2$. 
		
		Consider that the measures $\Lambda^i$ for $i=1,2$ are given (in polar coordinates) as
		\[
		\Lambda^i(A)= \int_0^{1}\int_{\mathbb{S}_+^2}1_{\{r\xi\in A\}}1_{\{r\xi\in\Delta\}}\lambda^i(d\xi)\frac{(1-\langle\xi,1\rangle u)^{\alpha-1}}{u^{1+\alpha}}du,
		\]
		with $\alpha\in(0,1)$. 
		Hence, we obtain that the process $Y$ is a solution to the martingale problem given in \eqref{mart_prob_iii_alt}, and therefore it 
		has the same distribution as the process $\overline{R}$ defined in Corollary \ref{cor_3}. Thus, the process $\overline{R}$ is a 2-type $\beta$-Fleming-Viot process with selection and mutation in the sense of \cite{Griff} (see for instance Example 2.3 in \cite{Griff} for the case with no mutation). This result clarifies the relation between $\alpha$-stable multi-type continuous-state branching processes and 2-type $\Lambda$-Fleming--Viot processes with selection and mutation, answering a question formulated in the concluding section of \cite{Griff}. By Theorem \ref{main} the time-change technique only works in the $\alpha$-stable case, however, we believe we can explain in general the relation between multi-type continuous-state branching processes and multi-type $\Lambda$-Fleming--Viot processes with selection and mutation, using the culling technique developed in \cite{CGP}. We leave it as a venue for future research.
	\end{remark}

\end{document}